\documentclass[12pt]{article}

\usepackage{times}
\usepackage{mathptm}
\usepackage{amsfonts}

\begin{document}

\title{Aftermath}
\author{Peter J. Cameron\\School of Mathematical Sciences\\Queen Mary,
University of London\\London E1 4NS, U.K.}
\date{}
\maketitle

\begin{abstract}
I take a quick overview at the recent development of combinatorics and its
current directions, as a discipline in its own right, as part of mathematics,
and as part of science and wider society.
\end{abstract}

Henry Whitehead reportedly said, ``Combinatorics is the slums of
topology''.\footnote{My attribution is confirmed by
Graham Higman, a student of Whitehead. Less disparagingly,
Hollingdale~\cite{holl} wrote ``\dots\ the branch of topology we now call
`graph theory' \dots''} This prejudice, the view that
combinatorics is quite different from `real mathematics', was not uncommon
in the twentieth century, among popular expositors as well as
professionals. In his biography of Srinivasa Ramanujan,
Robert Kanigel~\cite{kanigel} describes Percy MacMahon in these terms:
\begin{quote}\begin{raggedright}
[MacMahon's] expertise lay in combinatorics, a sort of 
glorified dice-throwing, and in it he had made contributions 
original enough to be named a Fellow of the Royal Society.
\end{raggedright}\end{quote}

In the later part of the century, attitudes changed.
When the 1998 film \textit{Good Will Hunting} featured a famous
mathematician at the Massachusetts Institute of Technology who
had won a Fields Medal for combinatorics, many found this
somewhat unbelievable.\footnote{The ``unsolvable math problem'' is based
on the actual experience of George B. Dantzig, who as a student solved two
problems posed by Jerzy Neyman at Berkeley in 1940: see
Brunvald~\cite{brunvald}.}
However, life followed art in this case
when, later in the same year, Fields Medals were awarded to Richard Borcherds
and Tim Gowers for work much of which was in combinatorics.

A more remarkable instance of life following art involves Stanis{\l}aw Lem's
1968 novel \textit{His Master's Voice}~\cite{lem}. The narrator, a
mathematician,  describes how he single-mindedly attacked his rival's work:
\begin{quote}
I do not think I ever finished any larger paper in all my younger work
without imagining Dill's eyes on the manuscript. What effort it cost me to
prove that the Dill variable combinatorics was only a rough approximation of
an ergodic theorem! Not before or since, I daresay, did I polish a thing so
carefully; and it is even possible that the whole concept of groups later
called Hogarth groups came out of that quiet, constant passion with which I
plowed Dill's axioms under.
\end{quote}
In 1975, Szemer\'edi~\cite{szem} published his remarkable combinatorial proof
that a set of natural numbers with positive density contains arbitrarily long
arithmetic progressions; in 1977, Furstenberg~\cite{furst} gave a proof based
on ergodic theory! (This is not to suggest that Furstenberg's attitude to
Szemer\'edi parallels Hogarth's to Dill in the novel.)

In this chapter, I have attempted to tease apart some of the
interrelated reasons for this change, and perhaps to throw some
light on present trends and future directions. I have divided
the causes into four groups:
the influence of the computer;
the growing sophistication of combinatorics;
its strengthening links with the rest of mathematics; and
wider changes in society.
I have told the story mostly through examples.

\section{The influence of the computer}

Even before computers were built, pioneers such as Babbage and Turing
realised that they would be designed on discrete principles, and would raise
theoretical issues which led to important mathematics. 

Kurt G\"odel~\cite{godel}
showed that there are true statements about the natural numbers
which cannot be deduced from the axioms of a standard system
such as Peano's. This result was highly significant for the
foundations of mathematics, but G\"odel's unprovable statement
itself had no mathematical significance. The first example of
a natural mathematical statement which is unprovable in Peano
arithmetic was discovered by Paris and Harrington~\cite{ph}, and
is a theorem in combinatorics (it is a slight strengthening of
Ramsey's theorem). It is unprovable from the axioms because the
corresponding `Paris--Harrington function' grows faster than
any provably computable function. Several further examples of
this phenomenon have been discovered, mostly combinatorial in
nature.\footnote{Calculating precise values for Ramsey
numbers, or even close estimates, appears to be one of the
most fiendishly difficult open combinatorial problems.}

More recently, attention has turned from \emph{computability} to
\emph{computational complexity}: given that something can be
computed, what resources (time, memory, etc.) are required for
the computation. A class of problems is said to be
\emph{polynomial-time computable}, or in $\mathsf{P}$, if any
instance can be solved in a number of steps bounded by a polynomial
in the input size. A class is in $\mathsf{NP}$ if the same
assertion holds if we are allowed to make a number of lucky
guesses (or, what amounts to the same thing, if a proposed solution
can be checked in a polynomial number of steps). The great unsolved
problem of complexity theory asks:
\begin{center}
Is $\mathsf{P}=\mathsf{NP}$?
\end{center}
On 24 May 2000, the Clay Mathematical Institute announced a list of
seven unsolved problems, for each of which a prize of one million dollars
was offered. The $\mathsf{P}=\mathsf{NP}$ problem was the first on the
list~\cite{clay}.

This problem is particularly important for combinatorics since
many intractable combinatorial problems (including the existence
of a Hamiltonian cycle in a graph) are known to be in
$\mathsf{NP}$. In the unlikely event of an affirmative solution,
`fast' algorithms would exist for all these problems.

Now we turn to the practical use of computers.

Computer systems such as \textsf{GAP}~\cite{gap} have been developed,
which can treat algebraic or combinatorial
objects, such as a group or a graph, in a way similar to the
handling of complex numbers or matrices in more traditional
systems. These give the mathematician a very powerful tool for
exploring structures and testing (or even formulating) conjectures.

But what has caught the public eye is the use of computers to
prove theorems. This was dramatically the case in 1976 when
Kenneth Appel and Wolfgang Haken~\cite{ah} announced that they
had proved the Four-Colour Theorem by computer. Their announcement
started a wide discussion over whether a computer proof is really
a `proof' at all: see, for example, Swart~\cite{swart} and Tymoczko~\cite{tym}
for contemporary responses. An even more massive computation by Clement
Lam and his co-workers~\cite{lametal}, discussed by Lam in~\cite{lam},
showed the non-existence of a projective plane of order~$10$.
Other recent achievements include the classification of Steiner triple
systems of order~$19$~\cite{ko}.

Computers have been used in other parts of mathematics. For
example, in the Classification of Finite Simple Groups (discussed
below), many of the sporadic simple groups were constructed with
the help of computers. The very practical study of fluid dynamics
depends on massive computation. What distinguishes combinatorics?
Two factors seem important:
\begin{itemize}
\item[(a)]  in a sense, the effort of the proof consists mainly in
detailed case analysis, or generate large amounts of data,
and so the computer does most of the work;
\item[(b)] the problem and solution are both discrete; the results
are not invalidated by rounding errors or chaotic behaviour.
\end{itemize}

Finally, the advent of computers has given rise to many new areas
of mathematics related to the processing and transmission of data.
Since computers are digital, these areas are naturally related
to combinatorics. They include coding theory (discussed below),
cryptography, integer programming, discrete optimisation, and
constraint satisfaction.

\section{The nature of the subject}

The last two centuries of mathematics have been dominated by the
trend towards axiomatisation. A structure which fails to satisfy the
axioms is not to be considered. (As one of my colleagues put it to
a student in a class, ``For a ring to pass the exam, it has to
get 100\%''.) Combinatorics has never fitted this pattern very well.

When Gian-Carlo Rota and various co-workers wrote an influential series of
papers with the title `On the foundations of combinatorial
theory' in the 1960s and 1970s (see~\cite{rota,cr}, for example),
one reviewer compared combinatorialists to nomads on the steppes
who had not managed to construct the cities in which other
mathematicians dwell, and expressed the hope that these papers
would at least found a thriving settlement.

While Rota's papers have been very influential, this view has not prevailed.
To see this, we turn to the more recent series on `Graph minors' by
Robertson and Seymour~\cite{rs}. These are devoted to the proof of a
single major theorem, that a minor-closed class of graphs is
determined by finitely many excluded minors. Along the way, a
rich tapestry is woven, which is descriptive (giving a topological
embedding of graphs) and algorithmic (showing that many graph
problems lie in $\mathsf{P}$) as well as deductive.

The work of Robertson and Seymour and its continuation is certainly
one of the major themes in graph theory at present, and has
contributed to a shorter proof of the Four-Colour Theorem, as well
as a proof of the Strong Perfect Graph Conjecture. Various authors,
notably Gerards, Geelen and Whittle, are extending it to classes of matroids
(see~\cite{ggw}).

What is clear, though, is that combinatorics will continue to
elude attempts at formal specification.

\section{Relations with mathematics}

In 1974, an Advanced Study Institute on Combinatorics was held at
Nijenrode, the Netherlands, organised by Marshall Hall and Jack van Lint.
This was one of the first presentations,
aimed at young researchers, of combinatorics as a mature mathematical
discipline. The subject was divided into five sections: theory of designs,
finite geometry, coding theory, graph theory, and combinatorial group theory.

It is very striking to look at the four papers in coding theory~\cite{hl}.
This was the youngest of the sections, having begun with the work
of Hamming and Golay in the late 1940s. Yet the methods being
used involved the most sophisticated mathematics: invariant theory,
harmonic analysis, Gauss sums, Diophantine equations.

This trend has continued. In the 1970s, the Russian school (notably Goppa,
Manin, and Vladut) developed links between coding theory and algebraic
geometry (specifically, divisors on algebraic curves). These links were
definitely `two-way', and both subjects benefited. More recently, codes
over rings and quantum codes have revitalised the subject and made new
connections with ring theory and group theory. In the related field of
cryptography, one of the most widely used ciphers is based on elliptic
curves.

Another example is provided by the most exciting development
in mathematics in the late 1980s, which grew 
from the work of Vaughan Jones, for which he received a 
Fields Medal in 1990. His research on traces of Von Neumann 
algebras came together with representations of the Artin
braid group to yield a new invariant of knots, with 
ramifications in mathematical physics and elsewhere. (See the
citation by Joan Birman~\cite{birman1} and her popular
account~\cite{birman2} for a map of this 
territory.) Later, it was pointed out that the Jones polynomial is
a specialisation of the Tutte polynomial, which had been defined
for arbitrary graphs by Tutte and Whitney and generalised to matroids
by Tutte. Tutte himself has given two accounts of his discovery:
\cite{tutte1,tutte2}. The connections led to further work. 
There was the work of Fran\c{c}ois Jaeger~\cite{jaeger}, who derived a 
spin model, and hence an evaluation of the Kauffman polynomial, 
from the strongly regular graph associated with the Higman--Sims 
simple group; and that of Dominic Welsh and his collaborators 
(described in his book~\cite{welsh}) on the computational complexity of 
the new knot invariants.

Sokal~\cite{sokal} has pointed out that there are close relations between
the Tutte polynomial and the partition function for the Potts model in
statistical mechanics; this interaction has led to important advances in
both areas.

Examples such as this of unexpected connections, by their nature,
cannot be predicted. However, combinatorics is likely to be involved
in such discoveries: it seems that deep links in mathematics often
reveal themselves in combinatorial patterns.

One of the best examples concerns the ubiquity of the Coxeter--Dynkin 
diagrams $A_n$, $D_n$, $E_6$, $E_7$, $E_8$. Arnol'd (see~\cite{arnold})
proposed finding an explanation of their ubiquity as a modern equivalent 
of a Hilbert problem, to guide the development of mathematics. 
He noted their occurrence in areas such as Lie algebras (the 
simple Lie algebras over $\mathbb{C}$), Euclidean geometry (root systems), 
group theory (Coxeter groups), representation theory (algebras of
finite representation type), and singularity theory (singularities 
with definite intersection form), as well as their connection 
with the regular polyhedra. To this list could be added 
mathematical physics (instantons) and combinatorics (graphs 
with least eigenvalue $-2$). Indeed, graph theory provides the 
most striking specification of the diagrams: they are just the
connected graphs with all eigenvalues smaller than $2$.

Recently this subject has been revived with the discovery by Fomin and
Zelevinsky~\cite{fz} of the role of the ADE diagrams in the theory of cluster
algebras: this is a new topic with combinatorial foundations and applications
in Poisson geometry, integrable systems, representation theory and total
positivity.

Other developments include the relationship of
combinatorics to finite group theory. The Classification of Finite
Simple Groups~\cite{gorenstein} is the greatest collaborative effort
ever in mathematics, running to about 15000 journal pages. (Ironically,
although the theorem was announced in 1980, the proof contained a gap
which has only just been filled.)
Combinatorial ideas (graphs, designs, codes, geometries) were involved
in the proof: perhaps most notably, the classification of spherical
buildings by Jacques Tits~\cite{tits}. Also, the result has had a
great impact in combinatorics, with consequences both for symmetric
objects such as  graphs and designs (see the survey by Praeger~\cite{praeger}),
and (more surprisingly) elsewhere as in Luks' proof~\cite{luks} that the graph
isomorphism problem for graphs of bounded valency is in $\mathsf{P}$.

This account would not be complete without a mention of the work of
Richard Borcherds~\cite{borcherds} on `monstrous moonshine', connecting
the Golay code, the Leech lattice, and the Monster simple group with
generalised Kac--Moody algebras and vertex operators in mathematical physics
and throwing up a number of product identities of the kind familiar from
the classic work of Jacobi and others.

\section{In science and in society}

Like any human endeavour, combinatorics has been affected by the great
changes in society last century. The first influence to be mentioned
is a single individual, Paul Erd\H{o}s, who is the subject of two
recent best-selling biographies~\cite{hoffman, schechter}. 

Erd\H{o}s' mathematical interests were wide, but combinatorics was
central to them. He spent a large part of his life without a permanent
abode, travelling the world and collaborating with hundreds of mathematicians.
In the days before email, he was a vital communication link between
mathematicians in the East and West; he also inspired a vast body of research
(his 1500 papers dwarf the output of any other modern mathematician).

Jerry Grossman~\cite{grossman} has demonstrated the growth in multi-author
mathematical papers this century, and how Erd\H{o}s was ahead of this
trend (and almost certainly contributed to it).

Erd\H{o}s also stimulated mathematics by publicising his vast collection of
problems; for many of them, he offered financial rewards for solutions. As an
example, here is one of his most valuable problems. Let $A=\{a_1,a_2,\ldots\}$
be a set of positive integers with the property that the sum of the
reciprocals of the members of $A$ diverges. Is it true that $A$ contains
arbitrarily long arithmetic progressions? The motivating special case
(recently solved affirmatively by Green and Tao~\cite{gt}) is that where $A$
is the set of prime numbers: this is a problem in number theory, but
Erd\H{o}s' extension to an arbitrary set transforms it into combinatorics.

Increased collaboration among mathematicians goes beyond the influence of
Erd\H{o}s; combinatorics seems to lead the trend. Aspects of this trend
include large international conferences (the Southeastern Conference on
Combinatorics, Graph Theory and Computing, which held its 42nd meeting in 2011,
attracts over 500 people annually), and electronic journals (the
\textit{Electronic Journal of Combinatorics}~\cite{ejc}, founded in 1994,
was one of the first refereed specialist electronic journals in mathematics).
Electronic publishing is particularly attractive to combinatorialists.
Often, arguments require long case analysis, which editors of traditional
print journals may be reluctant to include in full.

On a popular level, the Sudoku puzzle (a variant of the problem of completing
a critical set in a Latin square) engages many people in combinatorial
reasoning every day. Mathematicians have not been immune to its attractions.
At the time of writing, MathSciNet lists 38 publications with `Sudoku' in
the title, linking it to topics as diverse as spreads and reguli, neural 
networks, fractals, and Shannon entropy.

Our time has seen a change in the scientific
viewpoint from the continuous to the discrete. Two mathematical
developments of the twentieth century (catastrophe theory and
chaos theory) have shown how discrete effects can be produced by
continuous causes. (Perhaps their dramatic names reflect the
intellectual shock of this discovery.) But the trend is even
more widespread.

In their book introducing a new branch of discrete mathematics
(game theory), John von Neumann and Oskar Morgenstern~\cite{nm}
wrote:
\begin{quote}
The emphasis on mathematical methods seems to be shifted more
towards combinatorics and set theory -- and away from the
algorithm of differential equations which dominates mathematical
physics.
\end{quote}
How does discreteness arise in nature? Segerstr{\aa}le~\cite{us}
quotes John Maynard-Smith as saying
``today we really 
do have a mathematics for thinking about complex 
systems and things which undergo transformations 
from quantity into quality''
or from continuous to discrete, mentioning Hopf bifurcations as
a mechanism for this.

On the importance of discreteness in nature, Steven Pinker~\cite{pinker} 
has no doubt. He wrote:
\begin{quote}
It may not be a coincidence that the two systems in the universe
that most impress us with their open-ended complex design --
life and mind -- are based on discrete combinatorial systems.
\end{quote}
Here, `mind' refers primarily to language, whose combinatorial
structure is well described in Pinker's book. `Life' refers to
the genetic code, where DNA molecules can be regarded as words
in an alphabet of four letters (the bases adenine, cytosine, 
guanine and thymine), and three-letter subwords encode amino acids, 
the building blocks of proteins.

The Human Genome Project, whose completion was announced in 2001,
was a major scientific enterprise to
describe completely the genetic code of humans.
(See~\cite{bbklp} for an account of the mathematics involved, and 
\cite{lander2} for subsequent developments.)
At Pinker's university (the Massachusetts
Institute of Technology), the Whitehead Laboratory was
engaged in this project. Its director, Eric Lander, rounds off
this chapter and illustrates its themes. His doctoral thesis~\cite{lander}
was in combinatorics, involving a `modern' subject (coding theory),
links within combinatorics (codes and designs), and links to other
parts of mathematics (lattices and local fields). Furthermore,
he is a fourth-generation academic descendant of Henry Whitehead.

But there are now hints that discreteness plays an even more fundamental
role. One of the goals of physics at present is the construction of a theory
which could reconcile the two pillars of twentieth-century physics, general 
relativity and quantum mechanics. In describing string theory,
loop quantum gravity, and a variety of other approaches including 
non-commutative geometry and causal set theory, Smolin~\cite{smolin} argues
that all of them involve discreteness at a fundamental level (roughly the
Planck scale, which is much too small and fleeting to be directly observed).
Indeed, developments such as the holographic principle suggest that the basic
currency of the universe may not be space and time, but information,
measured in bits. Maybe the `theory of everything' will be combinatorial!

\end{document}